
\magnification=1200 %
\parindent 25pt %
\parskip=0pt %
\baselineskip=17pt %
\input amssym.def
\input amssym
\def\CC{\Bbb C}

\def\p2{${\bf P^2}$}
\def\sp{$H^0({\bf P^2}, {\cal O}(k))$}
 
\vbox to .92in { }
\centerline{\bf On the Intersections of Rational Curves with Cubic Plane Curves}
\vskip .2in
\centerline{Geng Xu %
\footnote*{ \sevenrm Partially Supported by NSF grant DMS-9596097}} 
\vskip .4in
\centerline{\bf 0. Introduction}
\par
\vskip .1in
{Let $V$ be a generic quintic threefold in the 4-dimensional complex projective space 
${\bf P^4}$. A well-known conjecture of Clemens says that $V$ has only 
a finite number of rational curves in each degree. Although Clemens'
conjecture is still quite open at this moment (it is known to be true for degree 
 up to 7 by S. Katz [K]), 
recently physicists, based on the mirror symmetry principle, have proposed a
formula to predict the numbers of rational curves of various degrees on $V$ which is a
Calabi-Yau threefold [M]. Furthermore,  
some of these predicted numbers have been verified
mathematically.}
\par
{Now if $S$ is a generic quartic surface in ${\bf P^3}$, then we know that 
for each positive integer $k$, there are  only a finite 
number of rational curves in $H^0(S, {\cal O}(k))$
 by the classification theory of surfaces. In this 
case,  the number of hyperplane sections with 3 nodes is known classically to be 3200 
(cf. [V], [YZ]).}
\par
{The main result of this paper is
 the following }
\par
\vskip .1in
\noindent
{\bf Theorem 1.} {\it Let $C$ be a smooth  curve on a smooth projective surface $S$ with
geometric genus $g(C) > 0$.
If $A$ is a line bundle on $S$ with $A\cdot C > 0$,
then there are only a finite number of reduced and irreducible
curves $D \in H^0(S,A)$ with geometric genus }
$$g(D) < {1\over 2} (K_S + C)A + 1$$
{\it such that
each $D$ intersects the  curve $C$ at exactly one point (set theoretically).  
Here the point where $D$ and $C$ intersect is not fixed.
}
\par
\vskip .1in
{Next we study one interesting example, we take $C$ to be a smooth cubic 
 on the projective plane \p2, which 
is a 1-dimensional ``Calabi-Yau'' manifold. 
In general, a degree $k$ plane curve intersects $C$ 
at $3k$ points by Bezout's Theorem. However, the set of plane rational curves of degree $k$ 
has dimension $3k-1$. Therefore, 
  one expects there are only a finite number of rational curves 
of degree $k$, each intersects $C$ at one distinct point for every positive integer $k$. We 
have} 
\par
\vskip .1in
\noindent
{\bf Corollary.} {\it If $C$ is a smooth cubic curve on the projective plane ${\bf P^2}$,  
then for every  positive integer $k$, there are only a finite number of rational 
curves of degree $k$, each intersects $C$ at exactly one point (set theoretically). 
Moreover, there exists such a rational curve for each $k$.}
\par
\vskip .1in
{The above corollary follows immediately from Theorem 1 and Proposition 2 
(see section 3). 
The number of such rational curves 
in the corollary is known when $k=1$, 
as a smooth cubic has 9 flexes. This number seems to 
be closely related to  the number of plane 
rational curves of degree $k$ passing through $3k-1$ general points (see  
section 3), which has been computed 
 (cf. [F1], [KM], [R1], 
 [RT]).} 
\par
{The method we use here is deformation of singularity. Namely, given a family of 
divisors in a linear system, we may get a new divisor by taking a derivative.}
\par
{Throughout this paper we work over the complex number field $\CC$. By a rational curve we 
mean a reduced and irreducible curve (possibly singular) with geometric genus 0.}  
\par
{Finally, I am very grateful to W. Cherry, I. Dogachev, M. Green,  
 Y. Kawamata and K. Oguiso for valuable discussions, and I would also like to 
thank the referee for helpful comments and suggestions.}
\par
\vskip .2in
\centerline{\bf 1. Singularity and deformation lemmas} 
\par
\vskip .1in
{In this section, we  first set up the necessary notations to 
describe general curve singularity, then we state a few lemmas which will 
be used in the next section.}
\par
{Let $D$ be a reduced and irreducible curve on a smooth projective surface 
$S$. If $P$ is a singular point of $D$, then there is a 
 desingularization of $D$ at $P$: }
$$S_{m+
1}(P) \buildrel \pi _{m+1} \over \longrightarrow S_m(P) \buildrel \pi _m
\over \longrightarrow \cdots \buildrel \pi _2 \over \longrightarrow S_1(P)
 \buildrel 
 \pi _1 \over \longrightarrow S_0(P) = S, $$
{so that the proper transform
$D^*_{m+1}$ of $D$ on $S_{m+1}(P)$ is smooth at every infinitely near point of 
$P$. Here $S_j(P) \buildrel \pi _j \over \longrightarrow 
 S_{j-1}(P) $ is the blow-up of $S_{j-1}(P)$ at a point $Q_{j-1}$ ($Q_0=P$)  
 with $E_{j-1}\subset S_j(P)$  the exceptional divisor. If 
 the proper transform $D^*_{j-1}$ of $D$   on $S_{
j-1}(P)$ has multiplicity $\mu_{j-1}$ at the point $Q_{j-1}$, that is, }
$$\pi^*_j (D^*_{j-1}) = D^*_j + \mu_{j-1} E_{j-1},$$
{then we say that $D$ has a {\it type $ \mu(P) = ( \mu _j, Q_j, E_j \mid j\in 
N(P)=\{
0,1,\ldots, m\})$ singularity at $P$}.}
\par
{Given a resolution of the singularity of $D\subset S$ at $P$ 
as above, if $F\subset
S$ is another curve (possibly reducible), such that}
$$ \pi ^*_j( \cdots (\pi ^*_2(\pi ^*_1(F) - \delta _0 E_0) - \delta _1 E_1) -  
\cdots ) - \delta _{j-1} E_{j-1} $$
{is an effective divisor for $j = 1, 2, \dots , m+1$ with some positive 
integers $\delta_i$, then we say that $F$ has a
{\it weak type $ \delta (P)= ( \delta _j, Q_j, E_j \mid j\in N(P))
$  singularity at $P$}. It is easy to see that a type $\mu (P)$ singularity 
 is a weak type $\mu (P)$ singularity at $P$.}
\par
{If $P_1, P_2, \cdots, P_r$ are all the distinct singular points of a 
curve $D$ on $S$, and $C$ has a type $\mu (P_i) = (\mu_{i,j}, Q_{i,j}, E_{i,j} 
| j\in N(P_i))$ singularity at $P_i$, then we say that $D$ has {\it a 
type $\mu  = (\mu_{i,j}, Q_{i,j}, E_{i,j} | i=1,\cdots,r, j\in N(P_i))$ 
singularity.} 
Likewise, if $F\subset S$ is another curve, it has a weak type 
$ \delta (P_i)= ( \delta _{i,j}, Q_{i,j}, E_{i,j} \mid j\in N(P_i))$ singularity 
at each $P_i$ $(i=1,\cdots,r)$, then we say that $F$ has {\it a weak type 
$\delta = (\delta_{i,j}, Q_{i,j}, E_{i,j} | i=1,\cdots,r, j\in N(P_i))$ singularity.}  
It is well-known that we have the following relationship  
between the geometric genus $g(D)$ and virtual genus $\pi(D)$ of $D$(cf. [I] ch. 9), }
$$\eqalignno{ g(D) &= \pi(D) - \sum_{i,j} {{ \mu_{i,j} (\mu_{i,j} 
- 1)}\over 2} \cr
&= {{ D\cdot (K_S + D)}\over 2} + 1 - \sum_{i,j} {{ \mu_{i,j} (\mu_{i,j}
- 1)}\over 2}. &(1)\cr}$$
\par
\vskip .1in
\noindent
{\bf Definition:} {\it Let $T\subset {\CC}$ be
 an open neighborhood of  the origin $0\in T$. Assume that $\sigma \colon D
\longrightarrow T$ is a family
of reduced and irreducible curves on a smooth projective surface $S$, 
$D_t=\sigma ^{-1}(t)$, then we say  
 that the family $D_t$ is equisingular at t=0 in the sense that we can  
resolve the
singularity of $D_t$ simultaneously, that is, there is a
proper morphism $\pi \colon \tilde D\longrightarrow D$,
 so that $\sigma \circ \pi \colon \tilde D \longrightarrow T$ is a flat map and
 }
$$\sigma \circ \pi \colon \tilde D_t = (\sigma \circ \pi )^{-1}(t) 
\longrightarrow 
 D_t$$
{\it is a resolution of the singularities of $D_t$. In addition, if $D_t$
has a type $\mu (t) = (\mu_{i,j}(t),$ $Q_{i,j}(t),E_{i,j}(t)\mid 
i=1,2,\cdots,r(t), j\in N(P_i(t))$ singularity
with the above resolution, then $\mu_{i,j}(t) = \mu_{i,j}$,  
$r(t)=r$ and $N(P_i(t))=N_i$  
are all independent of t, and the exceptional divisors and the singular loci of
the desingularization $\tilde D_t \longrightarrow D_t$ have the same 
configuration for all t (c.f. [Z],[W]).}
\par
\vskip .1in
{Let's now recall the following lemma from [X].}
\par
\vskip .1in
\noindent
{\bf Lemma 1.} {
\it Let $L$ be a line bundle on a smooth surface $S$, and $T\subset 
\CC$ be an open neighborhood of $0\in T$. Assume that $D_t = \{ G_t=0\} 
(t\in T)$ is a family of reduced and irreducible curves defined by global 
sections $G_t \in H^0(S, L)$ of $L$. If the family $D_t$ is equisingular 
at $t=0$, and the curve $D_t$ has a type $\mu (P(t)) = (\mu_j, Q_j(t), 
E_j(t) | j\in N) $ singularity at the point $P(t)\in S$, then the curve}
$$ \Bigl\{ {{dG_t}\over {dt}} \mid_{t=0} = 0\Bigr\} $$
{\it has a weak type $\mu (P(0))-1 = (\mu_j-1, Q_j(0), E_j(0) | j\in N) $  
 singularity at the point $P(0)\in S$.}  
\par
\vskip .1in
\noindent
{\it Proof.} {This is essentially Lemma 2.3 in [X]. $\quad \quad q. e. d.$}
\par
\vskip .1in
{The next lemma tells us how to compute the intersection number of 2 curves on a 
 smooth surface in terms of their singularities.}
\par
\vskip .1in
\noindent
{\bf Lemma 2.} 
{\it Assume that $D$ is a reduced and irreducible curve on a smooth surface $S$, 
$P_1, P_2, \cdots, P_r$ are all the distinct singular points of $D$ on $S$, and $D$ has a 
type $\mu  = (\mu_{i,j}, Q_{i,j}, E_{i,j} | i=1,\cdots,r, j\in N(P_i))$ 
singularity. Assume also that $F\subset S$ is another curve with $D\not\subset F$ which 
has  a weak type 
$\delta = (\delta_{i,j}, Q_{i,j}, E_{i,j} | i=1,\cdots,r, j\in N(P_i))$ singularity. Let}
$$S_{i,m(P_i)+
1} \buildrel \pi _{i,m(P_i)+1} \over \longrightarrow S_{i,m(P_i)} \buildrel \pi_{i,m(P_i)} 
\over \longrightarrow \cdots \buildrel \pi _{i,2} \over \longrightarrow S_{i,1}
 \buildrel 
 \pi _{i,1}     \over \longrightarrow S_{i,0} = S_{i-1,m(P_{i-1})+1}  $$
{\it be a desingularization of $D$ at $P_i$ $(i=1,2,\cdots,r)$ with 
 $S_{1,0}(P_1) = S$. 
Here we resolve the singularity of $D$ at $P_1$ first, then $P_2$, and so on. 
If we define effective divisors $F_0, F_1, \cdots, F_r$ inductively by}
$$\eqalign{ F_0 &= F, \cr
F_i &= \pi ^*_{i,m(P_i)+1}( \cdots (\pi ^*_{i,2}(\pi ^*_{i,1}(F_{i-1}) - \delta _{i,0} E_{i,0}) 
- \delta _{i,1}  E_{i,1}) -  
\cdots ) - \delta _{i,m(P_i)} E_{i,m(P_i)} \cr}$$
{\it for $i=1,2,\cdots,r$, then we have}
$$D\cdot F = \tilde D\cdot F_r + \sum_{i=1}^r \sum_{j=0}^{m(P_i)} \mu_{i,j} \delta_{i,j}. \eqno(2)$$
{\it Here $\tilde D$ is the proper transform of $D$ under the desingularization} 
 $$S_{r,m(P_r)+1} \longrightarrow \cdots \longrightarrow S_{1,0} = S.$$
\par
\vskip .1in
\noindent
{\it Proof.} 
{In general, if $\pi: X \longrightarrow S$ is the blow-up of $S$ at a point $x$ with 
exceptional divisor $E$, and $B\subset S$ is an irreducible curve and its multiplicity at $x$ 
is $mult_{x} B = e > 0$, then we have }
$$\pi^* B = B^* + e E, $$
{here $B^*$ is the proper 
transform of $B$ on $X$. Assume that $V$ is an effective divisor on $S$ so that 
$\pi^* V - \rho E$ is effective for some integer $\rho > 0$. Then we have }
$$\eqalign{B^* \cdot (\pi^* V - \rho E) &= (\pi^* B - e E)\cdot (\pi^* V - \rho E)\cr
&= B\cdot V - e \rho,\cr}$$
{that is,}
$$B\cdot V = B^* \cdot (\pi^* V - \rho E) + e \rho. \eqno(3)$$
\par
{Under the assumption of our lemma, the proper transform of $D$ on $S_{i,m(P_i)+1} $  
 and $D$ are isomorphic outside the points 
$P_1, \cdots, P_i$, therefore, they have the same type of 
singularity at the point $P_{i+1}$ $(i=1,2,\cdots,r-1)$. 
In a similar way, $F_i$ and $F$ have 
the same type of singularity at the point $P_{i+1}$. Thus,
 we deduce (2) by applying the above 
formula (3) repeatly to every step of the desingularization of $D$. $\quad \quad q. e. d.$} 
\par
\vskip .1in
{We end this section with the following elementary lemma.}
\par
\vskip .1in
\noindent
{\bf Lemma 3.} {\it Let $\Gamma $ be a $h$-dimensional linear system of 
divisors on a smooth projective variety with generators $\gamma_0, 
\gamma_1, \cdots, \gamma_h$. If $D_t = \sum b_i(t)\gamma_i \in \Gamma $ 
is a non-trivial family of divisors, then for generic $t$, the divisor}
$${d \over {dt}} D_t = \sum_{i=0}^h \bigl( {{db_i} \over {dt}} (t)\bigr)  
\gamma_i 
\not= \alpha (t) D_t,$$
{\it here $\alpha (t)$ is a function of $t$.}
\par
\vskip .1in 
\noindent
{\it Proof.} {If we view $\Gamma \cong {\bf P^h}$, then $D_t$ corresponds to 
a curve in ${\bf P^h}$. So, geometrically, Lemma 3 means that this curve has a 
tangent line at every smooth point. In detail, if }
$$
{d \over {dt}} D_t = \alpha (t) D_t$$
{for some functions $\alpha (t)$, 
then we have  
${{db_i} \over {dt}} (t) = \alpha (t) b_i(t),$
 hence}
$$b_i(t) = \beta_i e^{\int \alpha (t)dt}, \quad \quad i=0,1,\cdots,h$$
{for some constants $\beta_i$. Therefore, }
$$D_t = e^{\int \alpha (t)dt} (\sum_{i=0}^h \beta_i \gamma_i)$$ 
{are all equal to the same divisor $\sum \beta_i \gamma_i$. A 
contradiction. $\quad \quad q. e. d.$}

\par
\vskip .2in
\centerline{\bf 2. Proof of Theorem 1}
\par
\vskip .1in
{Let $C$ be a smooth  curve with geometric genus $g(C) > 0$, 
and $A$ be an ample line bundle on $C$ with deg $A = d > 0$.
 If $P\in C$ is a point 
 so that there is a divisor $D\in H^0(C, A)$  which is supported  
exactly at $P$, that is, $ D = d P$, 
then the following lemma, which is well-known,  tells us 
that the number of such points $P\in C$ is finite.}
\par
\vskip .1in
\noindent
{\bf Lemma 4.} {\it If $C $ is a smooth  curve with geometric genus $g(C) > 0$,  
and $A$ is an ample line bundle on $C$ with deg $A = d > 0$, then 
there are only a finite number of divisors 
$D\in H^0(C, A)$ having a point of multiplicity $d$.}
\par
\vskip .1in
{We now start the proof of Theorem 1.}
\par
\vskip .1in
\noindent
{\it Proof of Theorem 1.}  
{ Suppose that  there 
is a non-trivial family of   $D_t\in H^0(S, A)$  $(t\in T)$ with }
$$g(D_t)< {1\over 2} 
(K_S+C)A + 1$$
{such that each $D_t$ intersects $C$ at 
exactly one point, set theoretically. From Lemma 4 we know that, scheme theoretically,}
$$D_t\cap C = d R$$
{is independent of $t\in T$ for some point $R\in C$. Let $g\in H^0(S, A)$ be a global 
section for one of the $D_t$, and $\Gamma_R\subset  H^0(S, A)$ be the linear subspace generated 
by $g$ and $H^0(S, A-C)$. From the exact sequence}
$$0\longrightarrow H^0(S, A-C) \longrightarrow H^0(S, A) \longrightarrow H^0(C, A|_C),$$
{we know that $D_t\in \Gamma_R$ for all $t\in T$.}
\par
{As the virtual genus}
$$\pi (D_t) ={1\over 2} (K_S+D_t)D_t + 1 = {1\over 2}(K_S + A)A + 1$$
{is fixed, from (1) we know that the curves $D_t$ have only a finite type of
possible singularities. Without loss of generality, we may assume that the family $D_t$ is
equisingular at some generic point $t=t_0$,  $P_1(t), P_2(t), \cdots, P_r(t)$ are all the
distinct singular points of $D_t$ on $S$, and $D_t$ has a
type $\mu (t) = (\mu_{i,j},Q_{i,j}(t),E_{i,j}(t)\mid
i=1,2,\cdots,r, j\in N(P_i(t))$ singularity for $t$ in a neighborhood of $t_0$.
 }
\par
\vskip .1in
\noindent
{\bf Lemma 5.} 
{\it We have }  
$$A^2 \geq \sum_{i,j} \mu_{i,j}(\mu_{i,j} - 1) + A\cdot C.$$
\par
\vskip .1in
 {Assuming Lemma 5, we conclude from the genus formula (1) that }
$$\eqalign{ g(D_t) &= \pi (D_t) - \sum {{\mu_{i,j}(\mu_{i,j} -1)}\over 2}\cr
&\geq {1\over 2} (K_S+C)A + 1.\cr}$$
{A contradiction to our assumption $g(D_t)< {1\over 2} (K_S+C)A+1. \quad \quad q. e. d.$}  
\par
\vskip .1in
{Finally, we turn to the }
\par
\vskip .1in
\noindent
{\it Proof of Lemma 5.} {We prove Lemma 5 by induction on the number of steps to resolve 
the singularity of $D_t$ at $R$. For simplicity of notation, we assume that $t_0=0$.}
\par
\vskip .in

{\it 1.} 
 {The curve $D_t$ is smooth at $R$ for generic $t\in T$.  
 By Lemma 1, once we choose  some global sections  as  defining equations of $D_t$, 
}
$$F = {d\over {dt}}D_t|_{t=0} \in H^0(S, A)$$
{defines a curve $F$ on $S$, furthermore, $F$ has a weak type 
$\mu(0)-1 = (\mu_{i,j}-1, Q_{i,j}(0), E_{i,j}(0)  | $ $i=1,\cdots,r, j\in N(P_i))$ 
singularity. From 
Lemma 3, we have $D_0\not\subset F$, and $F\in \Gamma_R$ because 
$\Gamma_R$ is a linear space. }
\par
{If $C\not\subset F$, then we have that the intersection numbers }
$$I_R(C, D_0)=d,\quad  
I_R(C, F)\geq d. $$
{We conclude  that (cf. [F2] p. 7) }
$$I_R(D_0, F) \geq \hbox{min} (I_R(C, D_0), I_R(C, F)) \geq d.$$
{If $C\subset F$, we still have $I_R(D_0, F) \geq I_R(D_0, C)=d$.}  
\par
{Let $\tilde D_0 $ and $F_r$ be as in Lemma 2. Since $D_0$ is smooth at $R$ by assumption, 
$\tilde D_0$ and $D_0$ have the same local defining 
 equations near $R$, so do $F_r$ and $F$, we get}
$$\tilde D_0\cdot F_r \geq I_R(\tilde D_0, F_r) = I_R(D_0, F) \geq d.$$
{Consequently,}
$$\eqalign{A^2 = D_0\cdot F &= \tilde D_0\cdot F_r + \sum \mu_{i,j}(\mu_{i,j} -1) \cr
&\geq d + \sum \mu_{i,j}(\mu_{i,j} -1), \cr}$$
{thanks to (2). Therefore, the conclusion of Lemma 5 is true in this case. 
 We remark here that the above argument works also in the case $d=0$. }
\par
\vskip .1in
{\it 2.}
 {In general $R$ is a singular point of $D_t$ for generic $t\in T$. Assume that $P_1(t)=R$, and 
$P_2(t), \cdots, P_r(t)$ are the other singular points of $D_t$ on $S$. Now we continue to 
use the  notations  in step 1. We have}
$$I_R(C, D_t) = d \geq \mu_{1,0} = mult_R(D_t).$$
{Since $C$ is smooth at $R\in S$, we may choose local analytic coordinates $x,y$ near $R$ 
with $R=(0,0)$, so that $C$ is locally defined by $y=0$. Suppose that $D_t$ is locally 
defined by}
$$f_t(x, y) = x^{d} h_t(x) + y g_t(x, y),$$
{here the vanishing order of $g_t(x, y)$ at $R$ is $\mu_{1,0}-1$.  
 }
{Let}
$$S_{1,1} \buildrel \pi_{1,1} \over \longrightarrow S_{1,0} = S$$
{be the blow-up of $S$ at $R$ with exceptional divisor $E_{1,0}$ 
and new coordinates $x_1 = x, y_1=y/x$. 

Then the proper 
transform $D_{1,1}^*(t) \in H^0(S_{1,1}, \pi_{1,1}^*A-\mu_{1,0} E_{1,0})$
 of $D_t$ on $S_{1,1}$ is locally defined by}
$$x_1^{-\mu_{1,0}} f_t(x_1,x_1y_1) = x_1^{d-\mu_{1,0}} h_t(x_1) + y_1 (x_1^{-\mu_{1,0}+1} g_t
(x_1,x_1y_1)).\eqno(4)$$
{Let $R_{(1)}$ be the point on $S_{1,1}$ with coordinates $x_1=0, y_1=0$, and $C_{1,1}^*$ 
be the proper transform of $C$ on $S_{1,1}$. Then  $C_{1,1}^*$ is  
smooth at $R_{(1)}$ by assumption. From (4)  we have}
$$ C_{1,1}^*\cap D_{1,1}^*(t) = (d - \mu_{1,0})R_{(1)}. $$ 
\par
{Now the family $D_{1,1}^*(t)\in H^0(S_{1,1}, \pi_{1,1}^*A-\mu_{1,0} E_{1,0})$ has 
an improved 
singularity at the point $R_{(1)} \in S_{1,1}$, by induction we have}
$$(\pi_{1,1}^*A-\mu_{1,0} E_{1,0})^2 \geq \Bigl( \sum_{i=1}^r\sum_{j=0} 
\mu_{i,j}(\mu_{i,j} -1) -  \mu_{1,0}(\mu_{1,0}-1)\Bigr) + C_{1,1}^* \cdot 
(\pi_{1,1}^*A-\mu_{1,0} E_{1,0}).$$
{Meanwhile, }
$$\eqalign{ (\pi_{1,1}^*A-\mu_{1,0} E_{1,0})^2 &= A^2 - (\mu_{1,0})^2, \cr
C_{1,1}^* \cdot (\pi_{1,1}^*A-\mu_{1,0} E_{1,0}) &= (\pi_{1,1}^* C - E_{1,0})\cdot 
(\pi_{1,1}^*A-\mu_{1,0} E_{1,0})\cr
&= A\cdot C - \mu_{1,0}.\cr} $$
{Hence we deduce that  }
$$A^2\geq \sum_{i=1}^r \sum_{j=0} \mu_{i,j}(\mu_{i,j} -1) + A\cdot C.$$
{Therefore,  the conclusion of the lemma is also true for the family $D_t$.  
By induction, we are done. $\quad \quad q. e. d.$}
\par
\vskip .1in
\noindent
{\bf Remark 1.} { If $A$ is a line bundle on a smooth projective surface 
$S$, then by a similar argument to the  above, one shows that 
there are only a finite number of reduced and irreducible curves 
$D \in |A|$ with geometric genus $g(D) < {1\over 2} K_S \cdot A + 1$.}

\par
\vskip .2in
\centerline{\bf 3. The number of rational curves}
\par
\vskip .1in
{In this section, we take $C$ to be a smooth cubic curve on the projective plane $S=$\p2. 
First of all, it is well-known that the number of divisors in $H^0(C, {\cal O}_C(k))$
 having a point of multiplicity $3k$ is $9k^2$ (cf. [ACGH] p. 359).} 
\par
\vskip .1in
{In case $k=1$, we have that 
a smooth cubic has 9 ordinary flexes.}
\par
{When $k=2$,   
 the  number of divisors  in $H^0(C, {\cal O}(2))$ having a point of multiplicity 
6  is 36. Since a smooth conic is a rational curve, 
and $H^0({\bf P^2}, {\cal O}(2)) \cong H^0(C, {\cal O}(2))$,  
we have 27 smooth conics, each intersects $C$ at one point 
 because of the 9 double lines. 
}
\par
{When $k\geq 3$, let $P\in C$ be a point such that 
$(3k)P\in H^0(C, {\cal O}(k))$. Then the 
 number of such points $P\in C$ is $9k^2$. 
Fix such a point $P\in C$. Choose a 
homogeneous polynomial $g\in $\sp,  so that $\{ g=0\} \cap C = (3k)P$. Denote }
$$\Gamma_P^k = \{ s_1 g + s_2 g^* F | s_1, s_2, \in \CC, g^*\in H^0({\bf P^2}, {
\cal O}(k-3))
\} \subset H^0({\bf P^2}, {\cal O}(k)),$$
{here $F\in H^0({\bf P^2}, {\cal O}(3))$ is the defining equation of $C$. 
It is easy to see that $\Gamma_P^k $
 has
codimension 
 ${k+2\choose 2} - {k-3+2\choose 2} - 1 = 3k -1$
in \sp, and from the exact sequence }
$$0 \longrightarrow H^0({\bf P^2}, {\cal O}(k-3)) 
\longrightarrow H^0({\bf P^2}, {\cal O}(k)) \longrightarrow H^0(C, {\cal O}(k))
 \longrightarrow 0,$$ 
{we know that $D\in H^0({\bf P^2},{\cal O}(k))$ 
  has $D\cap C = (3k)P$ if and only if $D\in \Gamma_P^k$. }
\par
{Consider the Zariski closure $V_k$ of the set 
 of plane rational curves of degree $k$ on \p2 with ordinary  
nodes as the only singular points inside the space of all degree $k$ curves ${\bf P}H^0(
{\bf P^2}, {\cal O}(k))$. The degree of $V_k$, which can also be interpreted as the 
number of plane rational curves of degree $k$ passing through $3k-1$ general points, has 
been computed (cf. [F1], [KM], [R1], [RT]) to be, }
$$N_k = {1\over 2} \sum_{k_1+k_2=k} {{k_1k_2(3kk_1k_2-2k^2+6k_1k_2)(3k-4)!}\over 
{(3k_1-1)!(3k_2-1)!}} N_{k_1}N_{k_2}.$$ 
{Here $N_1=1, N_2=1, N_3=12, $ and $N_4=620$. 
 It is well-known that $V_k$ has dimension $3k-1$ in the space of all degree $k$ plane curves 
${\bf P}H^0({\bf P^2}, {\cal O}(k))$. Since $\Gamma_P^k$ has codimension $3k-1$, 
we have}
\par
\vskip .1in
\noindent
{\bf Proposition 1.} {\it The 
virtual number of divisors in $\Gamma_P^k\cap V_k$ is  
$N_k$, 
 with multiplicity suitably counted. }
\par
\vskip .1in
\par
{Of course, $\Gamma_P^k\cap V_k$ may contain reducible curves, and some curves may be counted 
more than once. There also might be rational curves with singular points worse than ordinary 
nodes which intersect the cubic curve $C$ at one point.  
Our next proposition is an easy consequence of the positivity of $N_k$. }
\par
\vskip .1in
\noindent
{\bf Proposition 2.} 
{\it If $C$ is a smooth cubic on the projective plane ${\bf P^2}$, then for every 
positive integer $k$, there exists a rational curve of degree $k$ which intersects 
$C$ at exactly one point.} 
\par
\vskip .1in
\noindent
{\it Proof.} 
{Let $J(C) = \CC /\Lambda $ be the Jacobian of $C$, $\omega \in H^0(C, \Omega^1)$
be a holomorphic differential and $p_0\in C$ be a fixed point. The map}
$$\varphi : C \longrightarrow J(C),\quad p \longrightarrow \int_{p_0}^p
\omega $$
{is an isomorphism. If we take any divisor $q_1+ q_2 + q_{3} \in H^0(C, {\cal O}(1))$,
by Abel's Theorem, for any point $P\in C$, $(3k) P\in H^0(C, {\cal O}(k))$ if and
only if }
$$(3k) \varphi (P) = \varphi ((3k)P) = \varphi (k(q_1 + q_2 + q_{3})) =
k(\varphi (q_1 + q_2 + q_{3})),$$
{that is,}
$$\varphi (P) \in \Bigl\{ {1\over {3}}  \varphi (q_1+q_2 +q_{3}) + {n\over {3k}}
\lambda_1 + {m\over {3k}} \lambda_2 \Bigr\}_{n,m\in \Bbb Z}, \eqno(5)$$
{here $\Lambda = \{ n\lambda_1 + m\lambda_2\}_{n,m\in \Bbb Z}$ is the lattice.}  
\par 
{For any positive integer $k$, if we choose    }
$$\varphi (P) =  {1\over {3}}  \varphi (q_1+q_2 +q_{3}) + {1\over {3k}}
\lambda_1 + {1\over {3k}} \lambda_2,$$
{then $(3k)P \in H^0(C, {\cal O}(k))$. From Proposition 1 we know  
there is a curve $D\in \Gamma_P^k \cap V_k$ because $N_k > 0$.
 Now we claim that $D$ is reduced and irreducible. 
 Otherwise, $D = D_1 D_2, $ deg$D_1 = k_1 < k, D_1 \cap C = (3k_1)P$, 
then}
$$\varphi (P) 
\in  \Bigl\{ {1\over {3}}  \varphi (q_1+q_2 +q_{3}) + {n_1\over {3k_1}}
\lambda_1 + {m_1\over {3k_1}} \lambda_2 \Bigr\}_{n_1,m_1\in \Bbb Z} $$
{by (5), which is impossible. }
\par
{Therefore, $D$ is a rational curve of degree $k$ which 
intersects $C$ at $P$. $\quad  q. e. d.$}  
\par
\vskip .1in
\noindent
{\it Question.} {How to calculate the exact number of rational curves of degree $k$,  
  each intersects a 
fixed smooth cubic $C$ at one point?}
\par
\vskip .1in
{We have an answer to the above question only in the first non-trivial case $k=3$.} 
\par 
{Let $C\subset $\p2 be a smooth cubic and $P$ is one of the 9 flexes of $C$. 
Then  one can choose homogeneous 
coordinates $X_0, X_1, X_2$ of \p2 so that $X_0 = 0$ is the inflection line 
tangent at $P=(0,0,1)$, and $C$ has the standard Weierstrass form (cf. [BK] p. 302)}
$$X_0X_2^2 = 4X_1^3 + \alpha X_0^2X_1 + \beta X_0^3.$$
{Consider the pencil}
$$\Gamma_P^3 = \{ s_1(4X_1^3 + \alpha X_0^2X_1 + \beta X_0^3 - X_0X_2^2) + s_2 X_0^3 
| s_1, s_2 \in \CC\} \subset H^0({\bf P^2}, {\cal O}(3)),$$
{generated by $C$ and the inflection line. 
 It is easy to see that  
$\Gamma_P^3$  contains 2 singular curves  
except when $\alpha =0$. If $\alpha = 0$, $C$ is isomorphic to the Fermat cubic 
$X_0^3 + X_1^3 + X_2^3 = 0$ because the $j-$invariants of both curves are 0, and 
$\Gamma_P^3$ contains only one singular curve which has a simple cusp.}
\par
\vskip .2in
{The study of the pencil $\Gamma_3(P)$ in the remaining case 
when $P$ is not a flex of the cubic $C$ was completed  by 
Keiji Oguiso after he saw the first version of this paper. 
The rest of this section is due to him.   }  
\par
{Consider a fixed point 
 $P\in C$ as in (5) and $P$ is not a flex of $C$. Since $C$ is smooth at $P$, there 
is a unique curve $D_P$ in the pencil $\Gamma_P^3$ such that $D_P$ is singular at $P$. 
By the genus formula,  $P$ must be a double point of $D_P$.  
One easily concludes that $P$ must be an ordinary double point 
(not a simple cusp) of $P$, otherwise $D_P$ will be reducible by local computations. }
\par
{From the proof of Proposition 2, we know that the pencil $\Gamma_P^3$ contains no 
reducible curves. 
Resolve the base locus of  $\Gamma_P^3$ by blowing up the projective plane 
\p2 9 times over the point $P$,  
}
$$\Psi : S\longrightarrow {\bf P^2}. $$
{Let $E_1, \cdots, E_8, E_9$ be the exceptional 
divisors of $\Psi $ in which the last one is $E_9$. Then $\Psi^* \Gamma_3(P) - 
\sum_{i=1}^9 E_i$ defines a relatively minimal elliptic fibration}
$$\phi : S \longrightarrow {\bf P^1}$$
{with a section $E_9$. Now the fibre of $\phi $ which contains the proper transform 
$D_P^* $ of $D_P$ is $D_P^* + E_1 + \cdots + E_8$, and the fact that $P$ is an 
ordinary double point of $D_P$ implies that $D_P^*$ intersects only $E_1$ and $E_8$. 
Since}
$$12 = \chi_{top} (S) = \sum \chi_{top} (\hbox{ singular fibers of }\phi ),$$
{and }
$$\chi (D_P^* + E_1 + E_2 + \cdots + E_8) = 9,$$
{the other singular fibers of $\chi $ must be either (a) 3 nodal rational curves, or (b) 
1 nodal rational curve and 1 cuspidal rational curve.}
\par
{Next, we exclude the case (b). 
Let $G_s = \phi^{-1}(s)-\{ sing  \phi^{-1}(s)\} $ for $s\in {\bf P^1}$. 
Then $G_s$ has a 
natural group structure (with origin $G_s\cap E_9$). If $\phi^{-1}(s_0) = C_1$ 
is a cuspidal rational curve, then 
$P$ is a smooth point of $C_1$ and we have $9P\equiv 3H \equiv 9T$ on $C_1$, here 
$T$ is the unique inflection point of $C_1$ and $H={\cal O}_{\bf P^2}(1)$. However,  
$G_s\cong G_a (\cong \CC )$, the additive group.  
 This implies that $P=T$ is the inflection point of $C_1$. Let $L_P$ be the inflection  
line of $C_1$ at $P$, since the intersection number  }
$$I_P(C, L_P) \geq \hbox{min } (I_P(C_1, L_P), 
I_P(C, C_1)) \geq 3,$$
{we conclude that $P$ is also an inflection point of $C$, which is 
impossible.}
\par
{In conclusion, $\Gamma_P^3$ contains 4 rational curves if $P$ is not a flex of $C$. 
Hence we have}
\par
\vskip .1in
\noindent
{\bf Proposition 3.} {\it The number of rational curves of degree 3 intersecting a 
smooth cubic $C$ at one point is either 306 (the general case) or 297. The number of 
such rational curves is 297 if and only if $C$ is projectively equivalent to  
 $X_0^3 + X_1^3 + X_2^3 = 0$ on \p2. }    
\par
\vskip .1in
\noindent
{\it Remark 2.} {After this work was done, we were informed that N. Takahashi [T] 
also calculated the number of degree 
3 rational curves in the generic case and he gave another proof of  
 our corollary in the introduction, 
 R. Vidunas [Vi] and Z. Ran [R2] obtained results similar to 
proposition 3.}
\par
\vskip .4in
 
\centerline{\bf References}

\par
\vskip .1in
\par
\noindent
\item{[ACGH]} E. Arbarello, M. Cornalba, P.A. Griffiths and J. Harris, 
{\it Geometry of Algebraic Curves}, Volume I, Springer-Verlag, 1985.
\par
\noindent
\item{[BK]} E. Brieskorn and H. Kn\"orrer, {\it Plane Algebraic Curves}, 
Birkh\"auser Verlag, 1986.  
\par
\noindent
\item{[F1]} W. Fulton, {\it Enumerative geometry via quantum cohomology; An 
introduction to the work of Kontsevich and Manin}, Summer Institute lecture notes 
(Santa Cruz), 1995. 
\par
\noindent
\item{[F2]} W. Fulton, {\it Intersection Theory}, Springer-Verlag, 1984.
\par
\noindent
\item{[GH]} P. Griffiths and J. Harris, {\it Principles of Algebraic Geometry},
John Wiley
and Sons, New York, 1978.

\par
\noindent
\item{[I]} S. Iitaka, {\it Algebraic Geometry}, Springer-Verlag, 1982.
\par
\noindent
\item{[K]} S. Katz, {\it On the finiteness of rational curves on quintic  
threefolds},
Compositio Math. {\bf 60} (1986), 151-162.

\par
\noindent
\item{[KM]} M. Kontsevich and Yu. Manin, {\it Gromov-Witten classes, 
quantum cohomology, and enumerative geometry}, Comm. Math. Phy. {\bf 164} 
(1994), 525-562. 
\par
\noindent
\item{[M]} D. Morrison, {\it Mirror symmetry and rational curves on quintic 
threefolds: a guide for mathematicians}, J. Amer. Math. Soc. {\bf 6} (1993), 
223-247.
\par
\noindent
\item{[R1]} Z. Ran, {\it Enumerative geometry of singular plane curves,} Inv. Math. 
{\bf 97} (1989), 447-465. 
\par
\noindent
\item{[R2]} Z. Ran, {\it The number of unisecant rational cubics to a 
plane cubic}, preprint.
\par
\noindent
\item{[RT]} Y. Ruan and G. Tian, {\it A mathematical theory of quantum 
cohomology}, J. Diff. Geom. {\bf 42} (1995), 259-367. 
\par
\noindent
\item{[T]} N. Takahashi, {\it Curves in the complement of a smooth plane 
cubic whose normalizations are $A^1$}, preprint.
\par
\noindent
\item{[V]} I. Vainsencher, {\it Enumeration of n-fold tangent hyperplanes to a 
surface}, J. Alg. Geom. {\bf 4} (1995),  503-526.
\par
\noindent
\item{[Vi]} R. Vidunas, {\it Rational curves intersecting an elliptic curve 
in ${\bf P^2}$ at one point and application to arithmetic geometry}, preprint.
\par
\noindent
\item{[W]} J. Wahl, {\it Equisingular deformations of plane algeboid curves},  
Trans. Amer. Math. Soc. {\bf 193} (1974), 143-170.

\par
\noindent
\item{[X]} G. Xu, {\it Subvarieties of general hypersurfaces in projective
space}, J. Diff. Geom. {\bf 39} (1994), 139-172.
\par
\noindent
\item{[YZ]} S.T. Yau and E. Zaslow, {\it BPS states, string duality, and 
nodal curves on K-3}, preprint.
\par
\noindent
\item{[Z]} O. Zariski, {\it Studies in equisingularity I. Equivalent 
singularities of plane algebroid curves,} Amer. J. Math. {\bf 87} (1965), 
507-536. 
\par
\vskip .2in
\noindent
{Department of Mathematics, Johns Hopkins University, Baltimore, MD 21218, USA}
\par
\noindent
E-mail: geng@math.jhu.edu

\end